\font\cmrfootnote=cmr10 scaled 750
\font\cmrhalf=cmr10 scaled \magstephalf

\font\cmrscfootnote=cmr7 scaled 750
\font\cmrschalf=cmr7 scaled \magstephalf

\font\cmrscscfootnote=cmr5 scaled 750
\font\cmrscschalf=cmr5 scaled \magstephalf

\font\mitfootnote=cmmi10 scaled 750
\font\mithalf=cmmi10 scaled \magstephalf

\font\mitscfootnote=cmmi7 scaled 750
\font\mitschalf=cmmi7 scaled \magstephalf

\font\mitscscfootnote=cmmi5 scaled 750
\font\mitscschalf=cmmi5 scaled \magstephalf

\font\cmsyfootnote=cmsy10 scaled 750
\font\cmsyhalf=cmsy10 scaled \magstephalf

\font\cmsyscfootnote=cmsy7 scaled 750
\font\cmsyschalf=cmsy7 scaled \magstephalf

\font\cmsyscscfootnote=cmsy5 scaled 750
\font\cmsyscschalf=cmsy5 scaled \magstephalf

\font\cmexfootnote=cmex10 scaled 750
\font\cmexhalf=cmex10 scaled \magstephalf

\font\cmexscfootnote=cmex10 scaled 750
\font\cmexschalf=cmex10 scaled \magstephalf

\font\cmexscscfootnote=cmex7 scaled 750
\font\cmexscschalf=cmex7

\def\mathfootnote{\textfont0=\cmrfootnote \textfont1=\mitfootnote \textfont2=\cmsyfootnote \textfont3=\cmexfootnote
        \scriptfont0=\cmrscfootnote \scriptscriptfont0=\cmrscscfootnote \scriptfont1=\mitscfootnote \scriptscriptfont1=\mitscscfootnote
        \scriptfont2=\cmsyscfootnote \scriptscriptfont2=\cmsyscscfootnote \scriptfont3=\cmexscfootnote \scriptscriptfont3=\cmexscscfootnote}
\def\mathhalf{\textfont0=\cmrhalf \textfont1=\mithalf \textfont2=\cmsyhalf \textfont3=\cmexhalf
        \scriptfont0=\cmrschalf \scriptscriptfont0=\cmrscschalf \scriptfont1=\mitschalf \scriptscriptfont1=\mitscschalf
        \scriptfont2=\cmsyschalf \scriptscriptfont2=\cmsyscschalf \scriptfont3=\cmexschalf \scriptscriptfont3=\cmexscschalf}

\font\Bbb=msbm10

\def\outlin#1{\hbox{\Bbb #1}}

\font\call=cmsy10 

\def\cal{\call}

\font\small=cmr5                       
\font\notsosmall=cmr7
\font\cmreight=cmr8

\font\cmrhalf=cmr10 scaled \magstephalf
       
\font\two=cmr10 scaled \magstep2

\font\sans=cmss10

\font\twosans=cmss10 scaled \magstep2

\font\caps=cmcsc10

\def\q{\quad}

\def\cut{\hfill\break}  \def\newline{\cut}
\def\h#1{\hbox{#1}}

\def\vp{\varphi}
  
\def\o{\omega}   \def\O{\Omega}

\let\H=\calH

\def\calO{{\h{\cal O}}} \def\calD{{\h{\cal D}}}

 \def\RR{{\outlin R}} \def\ZZ{{\outlin Z}} \def\NN{{\outlin N}}

\chardef\dotlessi="10  
\chardef\inodot="10

\def\polishL{\leavemode\setbox0=\hbox{L}\hbox to\wd0{\hss\char'40L}}

   \def\Ric{\h{\rm Ric}}

\def\GL2nR{\h{$GL(2n,\RR)$}}
\def\Sp2nR{\h{$Sp(2n,\RR)$}}

\def\part#1{\frac{\partial#1}{\partial t}}

\def\frac#1#2{{{#1}\over{#2}}}

\def\sm{\setminus}

\def\precpt
{\hbox{$\mskip3mu\mathhalf\subset\raise0.92pt\hbox{$\mskip-10mu\!\!\!\!
\mathfootnote\subset$}\mskip5mu$}}

\def\supsetnoteq{\hbox{$\mskip3mu\supset\raise-5.97pt
\hbox{$\mskip-10mu\!\!\!\scriptstyle\not=$}\mskip8mu$}}

\def\subsetnoteq{\hbox{$\mskip3mu\subset\raise-5.97pt
\hbox{$\mskip-11mu\!\!\scriptstyle\not=$}\mskip8mu$}}

\def\sseq{\subseteq}

\def\*{\star}

\def\D{\Delta}
 
\def\dbar{\bar\partial}
 \def\pa{\partial}
\def\ddbar{\partial\dbar}

\def\intM{\int_M}  \def\intm{\int_M}

\def\i{\sqrt{-1}}

\def\ra{\rightarrow}

\def\arrow#1{\hbox to #1pt{\rightarrowfill}}

\def\thhnotsosmall#1{${\hbox{#1}}^{\hbox{\small th}}$}

\def\MA{Monge-Amp\`ere }

\def\K{K\"ahler }
\def\Kno{K\"ahler}

\def\KE{K\"ahler-Einstein }  \def\KEno{K\"ahler-Einstein}

\def\po1{partition of unity }

\def\Loneloc{L_{\h{\small loc}}^1}

\def\strutdepth{\dp\strutbox}
\def\specialstar{\vtop to \strutdepth{
    \baselineskip\strutdepth
    \vss\llap{$\star$\ \ \ \ \ \ \ \ \  }\null}}
\def\marginalstar{\strut\vadjust{\kern-\strutdepth\specialstar}}
\def\marginal#1{\strut\vadjust{\kern-\strutdepth
    {\vtop to \strutdepth{
    \baselineskip\strutdepth
    \vss\llap{{ \small #1 }}\null} 
    }}
    }

\newcount\subsectionitemnumber
\def\clearsubsectionitemnumber{\subsectionitemnumber=0\relax}

\newcount\subsubsubsectionnumber
\def\clearsubsubsubsectionnumber{\subsubsubsectionnumber=0\relax}
\def\subsubsubsection#1{
\bigskip\noindent
\global\advance\subsubsubsectionnumber by 1%
{\rm
 \the\subsectionnumber.\the\subsubsectionnumber.\the\subsubsubsectionnumber}
{
#1.}
}

\newcount\subsubsectionnumber
\def\clearsubsubsectionnumber{\subsubsectionnumber=0\relax}
\def\subsubsection#1{
\clearsubsubsubsectionnumber
\bigskip\noindent
\global\advance\subsubsectionnumber by 1%
{%
\it \the\subsectionnumber.\the\subsubsectionnumber}
{
\it #1.}
}
\newcount\subsectionnumber
\def\clearsubsectionnumber{\subsectionnumber=0\relax}
\def\subsection#1{
\clearsubsectionitemnumber
\clearsubsubsectionnumber
\medskip\medskip\smallskip\noindent \global\advance\subsectionnumber by 1%
{%
\bf \the\subsectionnumber} 
{
\bf #1.}
}
\newcount\sectionnumber
\def\clearsectionnumber{\sectionnumber=0\relax}
\def\section#1{
\clearsubsectionnumber
\bigskip\bigskip\noindent \global\advance\sectionnumber by 1%
{%
\two
\the\sectionnumber} 
{
\two #1.}
}

\clearsectionnumber
\clearsubsectionnumber
\clearsubsubsectionnumber
\clearsubsubsubsectionnumber
\clearsubsectionitemnumber

\newcount\itemnumber
\def\clearitemnumber{\itemnumber=0\relax}
\def\c#1{ {\noindent\bf \the\itemnumber.} p. $#1$ \global\advance\itemnumber by 1}
\def\cn{ {\noindent\bf \the\itemnumber.} \global\advance\itemnumber by 1}

\clearitemnumber

\def\last#1{\advance\eqncount by -#1(\the\eqncount)\advance\eqncount by #1}
\def\llast{\advance\eqncount by -1(\the\eqncount)\advance\eqncount by 1}
\def\lllast{\advance\eqncount by -2(\the\eqncount)\advance\eqncount by 2}
\def\llllast{\advance\eqncount by -3(\the\eqncount)\advance\eqncount by 3}

\newcount\notenumber
\def\clearnotenumber{\notenumber=0\relax}
\def\note#1{\advance\notenumber by 1\footnote{${}^{(\the\notenumber)}$}
  {\lineskip0pt\notsosmall #1}}
\def\notewithcomma#1{\advance\notenumber by
1\footnote{${}^{\the\notenumber}$}
  {\lineskip0pt\notsosmall #1}}
                                                                                               
\clearnotenumber

\def\putnumber{%
\global\advance\subsectionitemnumber by 1{\the\subsectionnumber}.{\the\subsectionitemnumber}}
\def\numbering{{{\the\subsectionnumber}.{\the\subsectionitemnumber}}}

\def\Abstract#1{
{\narrower\bigskip\bigskip
\noindent {{\bf Abstract.\ \ } #1}

}
}

\def\Resume#1{
{\narrower\bigskip\bigskip
\noindent {{\bf R\'esum\'e.\ \ } #1}

}
}

\def\ref#1{{\bf[}{\sans #1}{\bf]}}
  % \ref with small skip

\def\sm{\smallskip}

\def\boxit#1{\vbox{\hrule\hbox{\vrule\kern3pt\vbox{\kern3pt#1\vglue3pt
\kern3pt}\kern3pt\vrule}\hrule}}

\long\def\frame#1#2#3#4{\hbox{\vbox{\hrule height#1pt
 \hbox{\vrule width#1pt\kern #2pt
 \vbox{\kern #2pt
 \vbox{\hsize #3\noindent #4}
\kern#2pt}
 \kern#2pt\vrule width #1pt}
 \hrule height0pt depth#1pt}}}

\def\help{\ifmmode\aftergroup\noindent\quad\else\quad\fi}
\def\helpp{\ifmmode\aftergroup\noindent\hfill\else\hfill\fi}

\def\on{\o^n}
\def\onminus1{\o^{n-1}/(n-1)!}
\def\ovp{\o_{\vp}}

\font\call=cmsy10   \font\Bbb=msbm10
\def\outlin#1{\hbox{\Bbb #1}}
\def\H{\hbox{\call H}}  
 
  \def\RR{{\outlin R}}
\def\frac#1#2{{{#1}\over{#2}}}
\newcount \eqncount
\def \eqnno{\global \advance \eqncount by 1 \futurelet \nexttok \parsenexttok}
\def \eqn{\global \advance \eqncount by 1 \eqno\futurelet \nexttok \parsenexttok}
\def \eqnd{\global\advance \eqncount by 1 \futurelet\nexttok\parsenexttokd}
\def \parsenexttok{\ifx \nexttok $\Nomark\else\expandafter \Mark\fi}
\def \parsenexttokd{\ifx \nexttok \hfil\Nomark\else\expandafter \Mark\fi}
\def \Nomark {(\the \eqncount)}\def \Mark #1{\xdef #1{(\the \eqncount)}#1}

\let\exa\expandafter
\catcode`\@=11 
\def\CrossWord#1#2#3{%
\def\@x@{}\def\@y@{#2}
\ifx\@x@\@y@ \def\@z@{#3}\else \def\@z@{#2}\fi
\exa\edef\csname cw@#1\endcsname{\@z@}}
\openin15=\jobname.ref
\ifeof15 \immediate\write16{No file \jobname.ref}%
   \else \input \jobname.ref \fi 
\closein15
\newwrite\refout
\openout\refout=\jobname.ref 
\def\warning#1{\immediate\write16{1.\the\inputlineno -- warning --#1}}
\def\Ref#1{%
\exa \ifx \csname cw@#1\endcsname \relax
\warning{\string\Ref\string{\string#1\string}?}%
    \hbox{$???$}%
\else \csname cw@#1\endcsname \fi}
\def\Tag#1#2{\begingroup
\edef\head{\string\CrossWord{#1}{#2}}%
\def\writeref{\write\refout}%
\exa \exa \exa
\writeref \exa{\head{\the\pageno}}%
\endgroup}

\catcode`\@=12

\def\Tagg#1#2{\Tag{#1}{#2}}

\def\TaggSection#1{\Tagg{#1}{Section~\the\subsectionnumber}}
\def\TaggS#1{\Tagg{#1}{\S~\the\subsectionnumber}}
\def\TaggSubS#1{\Tagg{#1}{\hbox{\S\S\unskip\the\subsectionnumber.\the\subsubsectionnumber}}}
\def\TaggSubsection#1{\Tagg{#1}{Subsection~\the\subsectionnumber.\the\subsubsectionnumber}}
\def\TaggEq#1{\Tagg{#1}{(\the\eqncount)}}
\def\Taggf#1{ \Tagg{#1}{{\bf[}{\sans #1}{\bf]}} }

\def\JJJ{\h{\rm J}}
\def\Ric{\hbox{\rm Ric}\,}

\def\Hc{\H_{c_1}}
\def\Dc{\calD_{c_1}}

\def\Hcplus{\H^{+}_{c_1}}

\def\O{\Omega}

\def\HO{\H_{\Omega}}

\def\Do{\calD_{\Omega}}

\def\MA{Monge-Amp\`ere }

\def\Vm{V^{-1}}

\magnification=1100
\hoffset1.1cm
\voffset1.3cm
\hsize5.09in
\vsize6.91in

\headline={\ifnum\pageno>1{\ifodd\pageno \oddheadline\else\evenheadline\fi}\fi}
\def\oddheadline{\centerline{\caps The Ricci iteration and its applications}}
\def\evenheadline{\centerline{\caps Y. A. Rubinstein}}

\overfullrule0pt
\parindent12pt
\vglue0.6cm

\noindent
\centerline{\twosans The Ricci iteration and its applications}
%\medskip
%\centerline{\twosans and its applications}
\smallskip

\font\tteight=cmtt8 

\bigskip
\centerline{\sans Yanir A. Rubinstein%
\footnote{$^*$}{\cmreight Massachusetts Institute of Technology. 
Email: {\tteight yanir@member.ams.org}
\hfill\break
\vglue-0.5cm
\hglue-\parindent\cmreight
Current address: Department of Mathematics, Princeton University, Princeton, NJ 08544.}}

\def\centereps#1#2#3{\vglue#2\relax\centerline{\hbox to#1%
            {\special{eps:#3.eps x=#1 y=#2}\hfil}}}

\font\cmreight=cmr8

\vglue0.03in

\footnote{}{\hglue-\parindent\cmreight April \thhnotsosmall{16}, 2007. Revised September 2007.
% Revised February \thhnotsosmall{28}, 2007.
\hfill\break 
\vglue-0.5cm
\hglue-\parindent\cmreight
\vglue-0.5cm 
\hglue-\parindent\cmreight
%       % 
58E11.} %  
%       % 
    
\vglue-0.12in

\bigskip
\Abstract{
In this Note we introduce and study dynamical systems related to the Ricci operator
on the space of K\"ahler metrics as discretizations of certain geometric 
flows. We pose a conjecture on their convergence towards canonical K\"ahler metrics
and study the case where the first
Chern class is negative, zero or positive. This construction has several applications in
K\"ahler geometry, among them an answer to a question of Nadel and a construction of multiplier
ideal sheaves. 
}

\Resume{
Dans cette Note nous introduisons et \'etudions des syst\`emes dynamiques reli\'ees 
\`a l'op\'erateur de Ricci sur l'espace des m\'etriques k\"ahl\'eriennes comme discr\'etisations
des certains flots g\'eom\'etriques. Nous posons une conjecture concernant leurs convergence vers
des m\'etriques k\"ahl\'eriennes canoniques and nous \'etudions le cas o\`u 
la premi\`ere classe de Chern est n\'egative, z\'ero ou positive. Cette construction a plusieurs
applications en g\'eom\'etrie k\"ahl\'erienne, parmi elles une r\'eponse \`a une question de Nadel et une
construction des faisceaux d'id\'eaux multiplicateurs.
}

\bigskip

\medskip
\subsectionnumber=0\relax

\subsection{Introduction}
\TaggSection{SectionIntroduction}
Our main purpose in this Note is to describe a new method for the construction
of canonical \K metrics via the discretization of certain geometric flows. The idea is
to turn a geometric flow into a set of difference equations.
Complete proofs will appear elsewhere \ref{12}.
 
The search for a canonical metric representative of 
a fixed \K class has been at the heart of \K geometry since its birth and has a long history
starting from \K and continuing, among many others, with the work of Calabi, Aubin, Yau, Tian
and Donaldson. For general extremal metrics a general existence theory is not presently available although
the so-called Yau-Tian-Donaldson conjecture suggests that it should be
related with notions of stability in algebraic geometry.

One of the main tools in the existence theory of \KE metrics and \Kno-Ricci solitons has been the Ricci flow introduced 
by Hamilton \ref{6}. Cao has shown that Yau's continuity method proof may be phrased in terms of the convergence of the Ricci flow \ref{3}. Much work has gone into understanding the Ricci flow on Fano
manifolds and recently Perelman and Tian and Zhu proved that the analogous convergence result holds in 
this case \ref{14}.
The idea that there might be another way of approaching canonical metrics, in the form of a discrete 
iterative dynamical system, 
was suggested by Nadel \ref{8} and is the main motivation for our work. 

\subsection{The Ricci iteration}
\TaggSection{SectionRicciIteration}
Let $(M,\JJJ,\o)$ be a connected compact closed \K manifold of complex dimension $n$
and let $\O$ denote a \K class. Let $\D_\o=-\dbar\circ\dbar^\star-\dbar^\star\circ\dbar$ denote
the Laplacian with respect to $\o$. Let $H_\o$ denote the Hodge projection operator from
the space of closed forms onto the kernel of $\D_\o$. Let $V=\O^n([M])$.
Denote by $\Do$ the space of all closed $(1,1)$-forms cohomologous to $\O$, by
$\HO$ the subspace of \K forms and by $\HO^+$ the 
subspace of \K forms whose Ricci curvature is positive.
We introduce the following
dynamical system on $\HO$ which is our main object of study in this Note.

\medskip
\noindent
{\bf Definition 2.1.}
{\it Given a \K form $\o\in\HO$ define the (time one) Ricci iteration%
\note{Most of the results hold also for discretizations corresponding to other time intervals.}
to be 
the sequence of metrics $\{\o_k\}_{k\ge0}$,
satisfying the following equations for each $k\in\NN$ for which a solution exists
$$
\o_{k+1}  =\o_k+H_{k+1}\Ric\o_{k+1}-\Ric\o_{k+1},\quad   k+1\in\NN,\quad
\o_0 =\o.
$$}
\medskip

This system of equations may be viewed as a discrete version of the flow
$\frac{\pa \o(t)}{\pa t}  =-\Ric\o(t)+H_t\Ric\o(t).$
This flow, first studied by Guan \ref{5}, can in turn be considered as a \K version of Hamilton's Ricci flow. 
Our work is motivated by the following conjecture (an analogue may also be posed for the flow):

\medskip
\noindent
{\bf Conjecture 2.2.}
{\it Let $(M,\JJJ)$ be a compact closed \K manifold, and assume that there
exists a constant scalar curvature \K metric in $\HO$.
Then for any $\o\in\HO$ the Ricci iteration exists for all $k\in\NN$ and converges
in an appropriate sense to a constant scalar curvature
metric.
}
\medskip

Our motivation for posing this conjecture comes from the following theorem.

\medskip

\noindent
{\bf Theorem 2.3. \it
Let $(M,\JJJ)$ be a compact closed \K manifold admitting a \KE metric. 
Let $\O$ be a \K class such that $\mu\O=c_1$ with $\mu\in\{0,\pm1\}$.
Then for any $\o\in\HO$ the Ricci iteration exists for all $k\in\NN$ and converges
in the sense of Cheeger-Gromov to a \KE metric.
}

\medskip

\noindent
{\bf Remark 2.4.} 
Note that the above conjecture may also be posed for solitons of the above flow, i.e., solutions of
$\hbox{\cal L}_X\o=\Ric\o-H_\o\Ric\o$ for a holomorphic vector field $X$, and
that a result analogous to Theorem 2.3 then holds for these metrics using the same methods and discretizing
the flow twisted by the one-parameter subgroup of automorphisms corresponding to $X$.  

\medskip

\noindent{\it Sketch of proof.}
First we prove that the iteration exists for each $k\in\NN$.
Since $H_\o\Ric\o = \mu \o$, this amounts to solving 
$
\o_{1} =\o_0+\mu\o_{1}-\Ric\o_{1}.
$
Let $\o_1=\o_{\vp_1}:=\o+\i\ddbar\vp_1$.
This can be written as a complex \MA equation:
$
\o_{\vp_1}^n=\o^n e^{f_\o-(\mu-1)\vp_1},\q \intm \o^n e^{f_\o-(\mu-1)\vp_1}=V,
$
where $\i\ddbar f_\o=\Ric\o-\o$, $\intm e^{f_\o}\on=V$.
The existence of solutions to such equations is known when $\mu\le0$ by the work of Aubin \ref{1} and Yau \ref{15}, 
and when $\mu=1$ by the work of Yau. Hence the iteration exists for each $k\in\NN$.

We divide the discussion into three cases, according to the sign of the first
Chern class. 

Assume first that $c_1<0$ and let $\O=-c_1$.
For each $k$ write $\o_k=\o_{\psi_k}$ with $\psi_k=\sum_{l=1}^k\vp_l$.
We have the following system of \MA equations:
$
\o_{\psi_k}^n
=
\o^n e^{f_\o+\psi_{k}+\vp_k},\,k\in \NN.
$
One readily sees that an inductive argument using the maximum principle implies that
$||\psi_{k}-\psi_{k-1}||_{C^{0}}\le C 2^{-k}$. This uniform bound implies
higher order a priori estimates, by elliptic regularity theory.
Therefore the sequence converges exponentially fast to a smooth function that 
we denote by $\psi_\infty$.

Consider the Chen-Tian functionals
$E_k$ \ref{4}.
We now observe the following monotonicity result. Its proof 
can be deduced from some of our previous results \ref{11}.

\medskip
\noindent
{\bf Lemma 2.5.} {\it Along the iteration $E_0$ is monotonically decreasing
whenever $\o_0\in\HO$. When $\mu=1$ the same is true for $E_1$, 
and if $\o_0\in\Hcplus$, also for $E_k,\, k\ge 2$.
}
\medskip

Coming back to the proof of the theorem, note that unless $\o_0$ is itself
\KEno, the functional $E_0$ is strictly decreasing along the iteration.
In particular since  $\o_\infty$ is a fixed point of the iteration it must be
\KEno.
The case $\mu=0$ is similar and so we omit the details.

Finally, we turn to the case $\mu=1$ and assume for simplicity 
that there are no holomorphic vector fields.
In the case $\mu=1$ the corresponding iteration takes a very special form.

\medskip

\noindent {\bf Definition 2.6.}
{\it
Define the
inverse Ricci operator $\Ric^{(-1)}:\Dc\ra\Hc$ by letting 
$\Ric^{(-1)}\o:=\o_\vp$ with $\ovp$ the unique \K form in $\Hc$ (given by the Calabi-Yau
Theorem \ref{15})
satisfying $\Ric\ovp=\o$. Similary denote higher order iterates of this operator
by $\Ric^{(-l)},\, l\in\ZZ$, with $\Ric^{(0)}:=\hbox{\rm Id}$.
} 

\medskip
 
We then see that the dynamical system for $\mu=1$
is nothing but the evolution of iterates of the inverse Ricci operator, $\o_k=\Ric^{(-k)}\o_0$.
For this case we are solving the system of equations
$
\o_{\psi_k}^n
=
\on e^{f_\o-\psi_{k-1}},\, k\in \NN.
$
Let $G_k$ be a Green function for $-\D_k:=-\D_{\o_{\vp_k}}$ 
satisfying $\intM G_k(\cdot,y)\o_{\psi_k}^n(y)=0$. Set $A_k=-\inf_{M\times M} G_k$.
Let $I(\o,\ovp)=\Vm\int_M\varphi(\o^n-\o^n_\vp)$.
Application of the Green formula gives
$
|\psi_k|\le n(A_0+A_k)+I(\o_0,\o_{\psi_k}).
$
Since $E_0$ is proper on $\Hc$ in the sense of Tian \ref{13}, 
we conclude that $I(\o,\o_{\psi_k})$ is uniformly bounded.
The crucial technical ingredient is now a uniform upper bound on the diameter.
Its derivation hinges on 
properties of the energy functionals, the definition of the iteration, and finally on an argument
due to Perelman adapted to this ``discrete" situation.
Now, combining the diameter estimate with the Green function estimate of Bando and Mabuchi \ref{2}, 
we conclude that $|\psi_k|\le C$. By monotonicity of 
the energy functionals we conclude that a subsequence can be chosen,  
converging to a \KE metric. This completes the 
outline of the proof of Theorem 2.3.

There are a number of applications of these constructions to several well-known objects of study
in \K geometry, among them canonical metrics, energy functionals, the Moser-Trudinger-Onofri
inequality, balanced metrics and the structure of the space of K\"ahler metrics.   
We point out two of the most obvious ones and describe others elsewhere. 
Also, these constructions can be generalized in some interesting directions \ref{12}.

The first application is an answer to a question raised by Nadel \ref{8}: Given 
$\o\in\Hc$ define 
a sequence of metrics $\o,\Ric\o,\Ric(\Ric\o),\ldots,$ as long as positivity 
is preserved; what are the periodic orbits of this dynamical system? The cases $k=2,3$ 
in the following theorem are due to Nadel.

\medskip

\noindent
{\bf Theorem 2.7.}
{\it
Let $(M,J,\o)$ be a Fano manifold and assume that $\Ric^{(k)}(\o)=\o$
for some $k\in\ZZ$. Then $\o$ is \KEno.
}

\medskip

\noindent
This is an immediate corollary of Lemma 2.5.

The second application is the construction of multiplier ideal sheaves \ref{7}
when a \KE metric does not exist. It may be seen as a discrete counterpart to a recent
result of Phong-\v Se\v sum-Sturm \ref{9}.

\medskip

\noindent
{\bf Theorem 2.8.}
{\it
Let $(M,\JJJ)$ be a Fano manifold that does not admit a \KE metric and let 
$\gamma>1$. One may extract a subsequence $\{\psi_{k_j}\}$ such that 
$\lim_{j\ra\infty}\psi_{k_j}=\psi_\infty$ exists in $\Loneloc(M)$ and
defines a nontrivial Nadel-type multiplier ideal sheaf
defined for each open set $U\sseq M$ by local sections
$\{h\in \calO_M(U): |h|^2 e^{-\gamma\psi_\infty}\in\Loneloc(M) \}$.
} 

\medskip

\vfill\eject
\def\smallblackbox{\vrule height.6ex width .5ex depth -.1ex}
\def\boxseparation{\hfil\smallblackbox$\q$\smallblackbox$\q$\smallblackbox\hfil}
\bigskip
\boxseparation
\bigskip

I would like to express my deep gratitude to my teacher, Gang Tian.
I thank X.-X. Chen, S. Donaldson, V. Guillemin, J. Morgan, T. Mrowka, I. Singer, J. Song and S. Zelditch for their interest in this work and for their warm encouragement.
Part of this research was carried out at Peking University in Summer 2005 where
the first version \ref{10} of \ref{12} was written and I thank that institution for its hospitality. 
This material is based upon work supported under a National Science 
Foundation Graduate Research Fellowship.

\frenchspacing

\bigskip\bigskip
\noindent{\bf Bibliography}
\bigskip
\def\ref#1{\Taggf{#1}\item{ {\bf[}{\sans #1}{\bf]} } }

\ref{1} {T. Aubin, \'{E}quations du type {M}onge-{A}mp\`ere sur les vari\'et\'es 
k\"ahl\'eriennes compactes, {\sl Bull. Sci. Math.} {\bf 102} (1978), 63--95.}
\sm

\ref{2} S. Bando, T. Mabuchi, Uniqueness of K\"ahler-Einstein metrics
modulo connected group actions, in: {\it Algebraic Geometry, Sendai, 1985,} 
Kinokuniya, 1987, 11--40.
\sm

\ref{3} H.-D. Cao, Deformations of K\"ahler metrics to K\"ahler-Einstein metrics on compact
K\"ahler manifolds, {\sl Inv. Math.} {\bf 81} (1985), 359--372.
\sm

\ref{4} X.-X. Chen, G. Tian, Ricci flow on K\"ahler-Einstein surfaces,
{\sl Inv. Math.} {\bf 147} (2002), 487--544.
\sm

\ref{5} D. Z.-D. Guan, Extremal-solitons and $C^\infty$ convergence of the modified Calabi flow on certain
$CP^1$ bundles, preprint, December 22, 2006.
\sm

\ref{6} R. S. Hamilton, Three-manifolds with positive Ricci curvature,
{\sl J. Diff. Geom.} {\bf 17} (1982), 255--306.
\sm

\ref{7} A. M. Nadel, Multiplier ideal sheaves and K\"ahler-Einstein metrics of positive scalar
curvature, {\sl Ann. Math.} {\bf 132} (1990), 549--596.
\sm

\ref{8} A. M. Nadel, On the absence of periodic points for the Ricci
curvature operator acting on the space of K\"ahler metrics,
in: {\it Modern Methods in Complex Analysis,} Princeton University Press, 1995, 277--281.
\sm

\ref{9} D. H. Phong, N. \v Se\v sum, J. Sturm,
Multiplier ideal sheaves and the K\"ahler-Ricci flow, preprint,
arxiv: math.DG/0611794.
\sm

\ref{10} Y. A. Rubinstein, On iteration of the Ricci operator on the space of K\"ahler metrics, I, manuscript,
August 14, 2005, unpublished.
\sm

\ref{11} Y. A. Rubinstein, On energy functionals, K\"ahler-Einstein metrics, and the Moser-Trudinger-Onofri 
neighborhood, preprint, arxiv: math.DG/0612440.
\sm

\ref{12} Y. A. Rubinstein, Some discretizations of geometric evolution equations and the Ricci iteration on the space of K\"ahler metrics, I, preprint, 2007; II, preprint (in preparation).
\sm

\ref{13} G. Tian,
K\"ahler-{E}instein metrics with positive scalar curvature, 
{\sl Inv. Math.} {\bf 130} (1997), {1--37}.
\sm

\ref{14} G. Tian, X.-H. Zhu, Convergence of K\"ahler-Ricci flow, 
{\sl J. Amer. Math. Soc.} {\bf 20} (2007), 675--699.
\sm

\ref{15} S.-T. Yau, On the Ricci curvature of a compact K\"ahler
manifold and the Complex Monge-Amp\'ere equation, I, {\sl Comm. Pure
Appl. Math.} {\bf 31} (1978), 339--411.

\end